\documentstyle[leqno]{article}
\newtheorem{th}{Theorem}[section]
\newtheorem{prop}[th]{Proposition}

\newcounter{defin}[section]
\renewcommand{\thedefin}{\thesection.\arabic{defin}}
\newcounter{ex}[section]
\renewcommand{\theex}{\thesection.\arabic{ex}}
\newcounter{rem}[section]
\renewcommand{\therem}{\thesection.\arabic{rem}}
\title{Upon $h$-normal $\Gamma$-linear connections\\
on $J^1(T,M)$}
\author{Mircea Neagu}
\date{}
\begin{document}
\maketitle
\begin{abstract}
Section 1 introduces the notion of $h$-normal $\Gamma$-linear connection on
the\linebreak 1-jet fibre bundle $J^1(T,M)$, and studies its local components. Section 2
analyses the main local components of torsion and curvature d-tensors attached
to an $h$-normal $\Gamma$-linear connection $\nabla$. Section 3 presents the
local Ricci identities  induced  by $\nabla$. The identities of the local
deflection d-tensors are also exposed. Section 4 is dedicated to the writing
of the local Bianchi identities of $\nabla$.
\end{abstract}

\section{Components of $h$-normal $\Gamma$-linear connections}

\hspace{5mm} Let $T$ (resp. $M$) be a {\it "temporal"} (resp. {\it
"spatial"}) manifold of dimension $p$ (resp. $n$) which is coordinated by
$(t^\alpha)_{\alpha=\overline{1,p}}$ (resp.$(x^i)_{i=\overline{1,n}}$). Let
us consider the 1-jet fibre bundle $J^1(T,M)\to T\times M$, naturally coordinated by
$(t^\alpha, x^i, x^i_\alpha)$. The  coordinate transformations on the product
manifold $T\times M$, induce the following coordinate transformations (gauge group)
on $J^1(T,M)$,
\begin{equation}
\left\{\begin{array}{l}
\tilde t^\alpha=\tilde t^\alpha(t^\beta)\\
\tilde x^i=\tilde x^i(x^j)\\
\tilde x^i_\alpha=\displaystyle{{\partial\tilde x^i\over\partial x^j}{\partial
t^\beta\over\partial\tilde t^\alpha}x^j_\beta}.
\end{array}\right.
\end{equation}
Note that, throughout this paper, the indices $\alpha,\beta,\gamma,\ldots$
run from $1$ to $p$ while the indices $i,j,k,\ldots$ run from $1$ to $n$.

On $E=J^1(T,M)$, we fixe a nonlinear connection $\Gamma$ defined
by the {\it temporal} components $M^{(i)}_{(\alpha)\beta}$ and the {\it spatial}
components $N^{(i)}_{(\alpha)j}$.We recall that the transformation rules of
the local components of the nonlinear connection $\Gamma$ are expressed
by \cite{10}
\begin{equation}
\left\{\begin{array}{l}\medskip
\displaystyle{\tilde M^{(j)}_{(\beta)\mu}{\partial\tilde t^\mu\over\partial
t^\alpha}=M^{(k)}_{(\gamma)\alpha}{\partial\tilde x^j\over\partial x^k}
{\partial t^\gamma\over\partial\tilde t^\beta}-{\partial\tilde x^j_\beta\over
\partial t^\alpha}}\\
\displaystyle{\tilde N^{(j)}_{(\beta)k}{\partial\tilde x^k\over\partial
x^i}=N^{(k)}_{(\gamma)i}{\partial\tilde x^j\over\partial x^k}
{\partial t^\gamma\over\partial\tilde t^\beta}-{\partial\tilde x^j_\beta\over
\partial x^i}}.
\end{array}\right.
\end{equation}
\addtocounter{ex}{1}
{\bf Example \theex} Let $h_{\alpha\beta}(t^\gamma)$ (resp. $g_{ij}(x^k)$) a
semi-Riemannian metric on the temporal (resp. spatial) manifold $T$ (resp. $M$),
and $H^\gamma_{\alpha\beta}$ (resp $\gamma^k_{ij}$) its Christoffel symbols.
Studying the transformation rules of the local components
\begin{equation}
\left\{\begin{array}{l}\medskip
M^{(j)}_{(\beta)\alpha}=-H^\gamma_{\alpha\beta}x^j_\gamma\\
N^{(j)}_{(\beta)i}=\gamma^j_{ik}x^k_\beta,
\end{array}\right.
\end{equation}
we conclude that $\Gamma_0=(M^{(j)}_{(\beta)\alpha},N^{(j)}_{(\beta)i})$
represents a nonlinear connection on $E$. This is called the {\it cannonical
nonlinear connection attached to the semi-Riemannian metrics $h_{\alpha\beta}$
and $\varphi_{ij}$}.\medskip

Let us consider $\displaystyle{\left\{{\delta\over\delta
t^\alpha}, {\delta\over\delta x^i}, {\partial\over\partial x^i_\alpha}\right\}
\subset{\cal X}(E)}$ and $\{dt^\alpha, dx^i, \delta x^i_\alpha\}\subset{\cal
X}^*(E)$ the adapted bases of the nonlinear connection $\Gamma$, where
\begin{equation}
\left\{\begin{array}{l}\medskip
\displaystyle{{\delta\over\delta t^\alpha}={\partial\over\partial t^\alpha}-
M^{(j)}_{(\beta)\alpha}{\partial\over\partial x^j_\beta}}\\\medskip
\displaystyle{{\delta\over\delta x^i}={\partial\over\partial x^i}-
N^{(j)}_{(\beta)i}{\partial\over\partial x^j_\beta}}\\
\delta x^i_\alpha=dx^i_\alpha+M^{(i)}_{(\alpha)\beta}dt^\beta+N^{(i)}_{(\alpha)
j}dx^j.
\end{array}\right.
\end{equation}
These  bases will be used in the description of geometrical objects on $E$,
because their transformation laws are very simple \cite{10}:
\begin{equation}\label{tlab}
\left\{\begin{array}{lll}\medskip
\displaystyle{{\delta\over\delta t^\alpha}={\partial\tilde t^\beta\over\partial
t^\alpha}{\delta\over\delta\tilde t^\beta}},&
\displaystyle{{\delta\over\delta x^i}={\partial\tilde x^j\over\partial
x^i}{\delta\over\delta\tilde x^j}},&
\displaystyle{{\partial\over\partial x^i_\alpha}={\partial\tilde x^j\over\partial
x^i}{\partial t^\alpha\over\partial\tilde t^\beta}{\partial\over\partial\tilde
x^j_\beta}},\\
\displaystyle{dt^\alpha={\partial t^\alpha\over\partial\tilde
t^\beta}d\tilde t^\beta},&
\displaystyle{dx^i={\partial x^i\over\partial\tilde
x^j}d\tilde x^j},&
\displaystyle{\delta x^i_\alpha={\partial x^i\over\partial\tilde
x^j}{\partial\tilde t^\beta\over\partial t^\alpha}\delta\tilde
x^j_\beta}.
\end{array}\right.
\end{equation}

In order to develope the theory of $\Gamma$-linear connections on the 1-jet
space $E$, we need the following
\begin{prop}
i) The Lie algebra ${\cal X}(E)$ of vector fields decomposes as
$$
{\cal X}(E)={\cal X}({\cal H}_T)\oplus{\cal X}({\cal H}_M)\oplus
{\cal X}({\cal V}),
$$
where
$$
{\cal X}({\cal H}_T)=Span\left\{{\delta\over\delta t^\alpha}\right\},\quad
{\cal X}({\cal H}_M)=Span\left\{{\delta\over\delta x^i}\right\},\quad
{\cal X}({\cal V})=Span\left\{{\partial\over\partial x^i_\alpha}\right\}.
$$

ii) The Lie algebra ${\cal X}^*(E)$ of covector fields decomposes as
$$
{\cal X}^*(E)={\cal X}^*({\cal H}_T)\oplus{\cal X}^*({\cal H}_M)\oplus
{\cal X}^*({\cal V}),
$$
where
$$
{\cal X}^*({\cal H}_T)=Span\{dt^\alpha\},\quad
{\cal X}^*({\cal H}_M)=Span\{dx^i\},\quad
{\cal X}^*({\cal V})=Span\{\delta x^i_\alpha\}.
$$
\end{prop}

Let us consider $h_T$, $h_M$ (horizontal) and $v$ (vertical) as the canonical projections of the above
decompositions. In this context, we have\medskip\\
\addtocounter{defin}{1}
{\bf Definition \thedefin} A linear connection $\nabla:{\cal X}(E)\times{\cal X}(E)\to
{\cal X}(E)$ is called {\it a $\Gamma$-linear connection on $E$} if $\nabla
h_T=0$, $\nabla h_M=0$ and $\nabla v=0$.\medskip

In order to describe in local terms a $\Gamma$-linear connection $\nabla$ on
$E$, we need nine unique local components,
\begin{equation}\label{lglc}\hspace*{5mm}
\nabla\Gamma=(\bar G^\alpha_{\beta\gamma},\;G^k_{i\gamma},\;G^{(i)(\beta)}_
{(\alpha)(j)\gamma},\;\bar L^\alpha_{\beta j},\;L^k_{ij},\;L^{(i)(\beta)}_{(\alpha)
(j)k},\;\bar C^{\alpha(\gamma)}_{\beta(k)},\;C^{j(\gamma)}_{i(k)},\;C^{(i)(\beta)
(\gamma)}_{(\alpha)(j)(k)}),
\end{equation}
which are locally defined by the relations
$$\begin{array}{l}\medskip
(h_T)\hspace{4mm}\displaystyle{
\nabla_{{\delta\over\delta t^\gamma}}{\delta\over\delta t^\beta}=
\bar G^\alpha_{\beta\gamma}{\delta\over\delta t^\alpha},\;\;
\nabla_{{\delta\over\delta t^\gamma}}{\delta\over\delta x^i}=G^k_{i\gamma}
{\delta\over\delta x^k},\;\;\nabla_{{\delta\over\delta t^\gamma}}
{\partial\over\partial x^i_\beta}=G^{(k)(\beta)}_{(\alpha)(i)\gamma}
{\partial\over\partial x^k_\alpha}},\\\medskip
(h_M)\hspace{3mm}\displaystyle{\nabla_{{\delta\over\delta x^j}}
{\delta\over\delta t^\beta}=\bar L^\alpha_{\beta j}{\delta\over\delta t^
\alpha},\;\;\nabla_{{\delta\over\delta x^j}}{\delta\over\delta x^i}=L^k_{ij}
{\delta\over\delta x^k},\;\;\nabla_{{\delta\over\delta x^j}}
{\partial\over\partial x^i_\beta}=L^{(k)(\beta)}_{(\alpha)(i)j}
{\partial\over\partial x^k_\alpha}},\\
(v)\hspace{6mm}\displaystyle{\nabla_{{\partial\over\partial x^j_\gamma}}
{\delta\over\delta t^\beta}=\bar C^{\alpha(\gamma)}_{\beta(j)}
{\delta\over\delta t^\alpha},\;\nabla_{{\partial\over\partial x^j_\gamma}}
{\delta\over\delta x^i}=C^{k(\gamma)}_{i(j)}
{\delta\over\delta x^k},\;\nabla_{{\partial\over\partial x^j_\gamma}}
{\partial\over\partial x^i_\beta}=C^{(k)(\beta)(\gamma)}_{(\alpha)(i)(j)}
{\partial\over\partial x^k_\alpha}}.
\end{array}
$$
\addtocounter{rem}{1}
{\bf Remark \therem} The transformation rules of the above connection coefficients
are completely described in \cite{12}.\medskip\\
\addtocounter{ex}{1}
{\bf Example \theex} Let us consider $\Gamma_0=(M^{(i)}_{(\alpha)\beta},
N^{(i)}_{(\alpha)j})$ the canonical nonlinear connection on $E$, attached to
the semi-Riemannian metrics pair $(h_{\alpha\beta},\varphi_{ij})$. In these
conditions, the following local coefficients \cite{12}
$$
B\Gamma_0=(\bar G^\alpha_{\beta\gamma},0,G^{(k)(\beta)}_{(\alpha)(i)\gamma},
0,L^k_{ij},L^{(k)(\beta)}_{(\alpha)(i)j},0,0,0),
$$
where $\bar G^\gamma_{\alpha\beta}=H^\gamma_{\alpha\beta},\;G^{(k)(\beta)}_
{(\gamma)(i)\alpha}=-\delta^k_iH^\beta_{\alpha\gamma},\;L^k_{ij}=\gamma_{ij}
^k$ and $L^{(k)(\beta)}_{(\gamma)(i)j}=\delta^\beta_\gamma\gamma^k_{ij}$,
verify the transformation rules of the local coefficients of a $\Gamma_0$-linear
connection. This is called {\it the Berwald
$\Gamma_0$-linear connection of the metrics pair $(h_{\alpha\beta},\varphi_
{ij})$}.\medskip

Now, let $\nabla$ be a $\Gamma$-linear connection on $E$, locally defined by
\ref{lglc}.
The linear connection $\nabla$ induces a natural linear connection on the
d-tensors set of the jet fibre bundle $E=J^1(T,M)$, in the following fashion:
starting with a vector field $X$ and a d-tensor field $D$
locally expressed by
$$\begin{array}{l}\medskip\displaystyle{
X=X^\alpha{\delta\over\delta t^\alpha}+X^m{\delta\over\delta x^m}
+X^{(m)}_{(\alpha)}{\partial\over\partial x^m_\alpha},}\\
\displaystyle{
D=D^{\alpha i(j)(\delta)\ldots}_{\gamma k(\beta)(l)\ldots}{\delta\over\delta
t^\alpha}\otimes{\delta\over\delta x^i}\otimes{\partial\over\partial x^j_
\beta}\otimes dt^\gamma\otimes dx^k\otimes\delta x^l_\delta\ldots,}
\end{array}
$$
we introduce the covariant derivative
$$
\begin{array}{l}\medskip
\nabla_XD=X^\varepsilon\nabla_{\delta\over\delta t^\varepsilon}D+X^p\nabla_
{\delta\over\delta x^p}D+X^{(p)}_{(\varepsilon)}\nabla_{\partial\over\partial
x^p_\varepsilon}D=\left\{X^\varepsilon D^{\alpha i(j)(\delta)\ldots}_{\gamma k
(\beta)(l)\ldots /\varepsilon}+X^p\right.\\\left.
D^{\alpha i(j)(\delta)\ldots}_{\gamma k(\beta)(l)\ldots\vert p}
+X^{(p)}_{(\varepsilon)}D^{\alpha i(j)(\delta)\ldots}_{\gamma k (\beta)(l)
\ldots}\vert_{(p)}^{(\varepsilon)}\right\}\displaystyle{
{\delta\over\delta t^\alpha}\otimes{\delta\over\delta x^i}\otimes{\partial
\over\partial x^j_\beta}\otimes dt^\gamma\otimes dx^k\otimes\delta x^l_\delta
\ldots,}
\end{array}
$$
where\\\\
$
(h_T)\hspace{6mm}\left\{\begin{array}{l}\medskip\displaystyle{
D^{\alpha i(j)(\delta)\ldots}_{\gamma k(\beta)(l)\ldots /\varepsilon}=
{\delta D^{\alpha i(j)(\delta)\ldots}_{\gamma k(\beta)(l)\ldots}\over
\delta t^\varepsilon}+D^{\mu i(j)(\delta)\ldots}_{\gamma k(\beta)(l)\ldots}
\bar G^\alpha_{\mu\varepsilon}+}\\\medskip
D^{\alpha m(j)(\delta)\ldots}_{\gamma k(\beta)(l)\ldots}G^i_{m\varepsilon}+
D^{\alpha i(m)(\delta)\ldots}_{\gamma k(\mu)(l)\ldots}G^{(j)(\mu)}_{(\beta)
(m)\varepsilon}+\ldots-\\
-D^{\alpha i(j)(\delta)\ldots}_{\mu k(\beta)(l)\ldots}\bar G^\mu_{\gamma
\varepsilon}-D^{\alpha i(j)(\delta)\ldots}_{\gamma m(\beta)(l)\ldots} G^m_
{k\varepsilon}-D^{\alpha i(j)(\mu)\ldots}_{\gamma k(\beta)(m)\ldots} G^{(m)
(\delta)}_{(\mu)(l)\varepsilon}-\ldots,
\end{array}\right.\medskip
$
$
(h_M)\hspace{5mm}\left\{\begin{array}{l}\medskip\displaystyle{
D^{\alpha i(j)(\delta)\ldots}_{\gamma k(\beta)(l)\ldots\vert p}=
{\delta D^{\alpha i(j)(\delta)\ldots}_{\gamma k(\beta)(l)\ldots}\over
\delta x^p}+D^{\mu i(j)(\delta)\ldots}_{\gamma k(\beta)(l)\ldots}
\bar L^\alpha_{\mu p}+}\\\medskip
D^{\alpha m(j)(\delta)\ldots}_{\gamma k(\beta)(l)\ldots}L^i_{mp}+
D^{\alpha i(m)(\delta)\ldots}_{\gamma k(\mu)(l)\ldots}L^{(j)(\mu)}_{(\beta)
(m)p}+\ldots-\\
-D^{\alpha i(j)(\delta)\ldots}_{\mu k(\beta)(l)\ldots}\bar L^\mu_{\gamma
p}-D^{\alpha i(j)(\delta)\ldots}_{\gamma m(\beta)(l)\ldots} L^m_
{kp}-D^{\alpha i(j)(\mu)\ldots}_{\gamma k(\beta)(m)\ldots} L^{(m)
(\delta)}_{(\mu)(l)p}-\ldots,
\end{array}\right.\medskip
$
$
(v)\hspace{8mm}\left\{\begin{array}{l}\medskip\displaystyle{
D^{\alpha i(j)(\delta)\ldots}_{\gamma k(\beta)(l)\ldots}\vert_{(p)}^{(\varepsilon)}=
{\partial D^{\alpha i(j)(\delta)\ldots}_{\gamma k(\beta)(l)\ldots}\over
\partial x^p_\varepsilon}+D^{\mu i(j)(\delta)\ldots}_{\gamma k(\beta)(l)\ldots}
\bar C^{\alpha(\varepsilon)}_{\mu(p)}+}\\\medskip
D^{\alpha m(j)(\delta)\ldots}_{\gamma k(\beta)(l)\ldots}C^{i(\varepsilon)}_{m(p)}+
D^{\alpha i(m)(\delta)\ldots}_{\gamma k(\mu)(l)\ldots}C^{(j)(\mu)(\varepsilon)}
_{(\beta)(m)(p)}+\ldots-\\
-D^{\alpha i(j)(\delta)\ldots}_{\mu k(\beta)(l)\ldots}\bar C^{\mu(\varepsilon)}_
{\gamma(p)}-D^{\alpha i(j)(\delta)\ldots}_{\gamma m(\beta)(l)\ldots} C^{m(\varepsilon)}
_{k(p)}-D^{\alpha i(j)(\mu)\ldots}_{\gamma k(\beta)(m)\ldots} C^{(m)(\delta)
(\varepsilon)}_{(\mu)(l)(p)}-\ldots.
\end{array}\right.\medskip
$

The local operators "$ _{/\varepsilon}$", "$_{\vert p}$" and "$\vert
^{(\varepsilon)}_{(p)}$" are called the {\it $T$-horizontal covariant derivative,
$M$-horizontal covariant derivative} and {\it vertical covariant derivative}
of the $\Gamma$-linear connection $\nabla$.\medskip\\
\addtocounter{rem}{1}
{\bf Remarks \therem} i) In the particular case of a function $f(t^\gamma,x^k,x^k_
\gamma)$ on $J^1(T,M)$, the above covariant derivatives reduce to
\begin{equation}
\left\{\begin{array}{l}\medskip
f_{/\varepsilon}=\displaystyle{{\delta f\over\delta t^\varepsilon}={\partial f
\over\partial t^\varepsilon}-M^{(k)}_{(\gamma)\varepsilon}{\partial f\over
\partial x^k_\gamma}}\\\medskip
f_{\vert p}=\displaystyle{{\delta f\over\delta x^p}={\partial f
\over\partial x^p}-N^{(k)}_{(\gamma)p}{\partial f\over\partial x^k_\gamma}}\\
f\vert^{(\varepsilon)}_{(p)}=\displaystyle{{\partial f\over\partial x^p_
\varepsilon}}.
\end{array}\right.
\end{equation}

ii) Particularly, starting with a d-vector field $X$ on $J^1(T,M)$, locally
expressed by
$$
X=X^\alpha{\delta\over\delta t^\alpha}+X^i{\delta\over\delta x^i}+
X^{(i)}_{(\alpha)}{\partial\over\partial x^i_\alpha},
$$
the following expressions of above covariant derivatives hold good:\medskip\\
$(h_T)\hspace*{6mm}
\left\{\begin{array}{l}\medskip
\displaystyle{X^\alpha_{/\varepsilon}={\delta X^\alpha\over\delta t^\varepsilon}
+X^\mu\bar G^\alpha_{\mu\varepsilon}}\\\medskip
\displaystyle{X^i_{/\varepsilon}={\delta X^i\over\delta t^\varepsilon}
+X^mG^i_{m\varepsilon}}\\
\displaystyle{X^{(i)}_{(\alpha)/\varepsilon}={\delta X^{(i)}_{(\alpha)}\over
\delta t^\varepsilon}+X^{(m)}_{(\mu)}G^{(i)(\mu)}_{(\alpha)(m)\varepsilon}},
\end{array}\right.\medskip
$
$(h_M)\hspace*{5mm}
\left\{\begin{array}{l}\medskip
\displaystyle{X^\alpha_{\vert p}={\delta X^\alpha\over\delta x^p}
+X^\mu\bar L^\alpha_{\mu p}}\\\medskip
\displaystyle{X^i_{\vert p}={\delta X^i\over\delta x^p}
+X^mL^i_{mp}}\\
\displaystyle{X^{(i)}_{(\alpha)\vert p}={\delta X^{(i)}_{(\alpha)}\over
\delta x^p}+X^{(m)}_{(\mu)}L^{(i)(\mu)}_{(\alpha)(m)p}},
\end{array}\right.\medskip
$
$(v)\hspace*{8mm}
\left\{\begin{array}{l}\medskip
\displaystyle{X^\alpha\vert^{(\varepsilon)}_{(p)}={\partial X^\alpha\over
\partial x^p_\varepsilon}+X^\mu\bar C^{\alpha(\varepsilon)}_{\mu(p)}}\\\medskip
\displaystyle{X^i\vert^{(\varepsilon)}_{(p)}={\partial X^i\over\partial x^p_
\varepsilon}+X^mC^{i(\varepsilon)}_{m(p)}}\\
\displaystyle{X^{(i)}_{(\alpha)}\vert^{(\varepsilon)}_{(p)}={\partial X^{(i)}_
{(\alpha)}\over\partial x^p_\varepsilon}+X^{(m)}_{(\mu)}C^{(i)(\mu)(\varepsilon)}
_{(\alpha)(m)(p)}}.
\end{array}\right.
$

iii) The local covariant derivatives associated to the Berwald $\Gamma_0$-linear
connection, will be denoted by
"$ _{/\!/\varepsilon}$", "$_{\Vert p}$" and "$\Vert^{(\varepsilon)}_{(p)}$".
\medskip

Now, let $h_{\alpha\beta}$ be a fixed pseudo-Riemannian metric on the temporal
manifold $T$, $H^\gamma_{\alpha\beta}$ its Christoffel symbols and $J=J^{(i)}
_{(\alpha)\beta j}{\partial\over\partial x^i_\alpha}\otimes dt^\beta\otimes
dx^j$, where $J^{(i)}_{(\alpha)\beta j}=h_{\alpha\beta}\delta^i_j$, the
{\it normalization d-tensor} \cite{10} attached to the metric $h_{\alpha\beta}$.
The big number of coefficients which characterize a $\Gamma$-linear
connection on $E$, determines us to consider the following\medskip\\
\addtocounter{defin}{1}
{\bf Definition \thedefin} A $\Gamma$-linear connection $\nabla$ on $J^1(T,M)$,
defined by the local coefficients
$$
\nabla\Gamma=(\bar G^\alpha_{\beta\gamma},G^k_{i\gamma},G^{(k)(\beta)}_
{(\alpha)(i)\gamma},\bar L^\alpha_{\beta j},L^k_{ij},L^{(k)(\beta)}_{(\alpha)
(i)j},\bar C^{\alpha(\gamma)}_{\beta(j)},C^{k(\gamma)}_{i(j)},C^{(k)(\beta)
(\gamma)}_{(\alpha)(i)(j)}),
$$
that verify the relations $\bar G^\alpha_{\beta\gamma}=H^\alpha_{\beta\gamma},
\;\bar L^\alpha_{\beta j}=0,\;\bar C^{\alpha(\gamma)}_{\beta(j)}=0$ and
$\nabla J=0$, is called an {\it $h$-normal $\Gamma$-linear connection}.\medskip\\
\addtocounter{rem}{1}
{\bf Remark \therem} Taking into account the local covariant $T$-horizontal "$_{
/\gamma}$", $M$-horizontal "$_{\vert k}$" and vertical "$\vert^{(
\gamma)}_{(k)}$" covariant derivatives induced by $\nabla$, the condition
$\nabla J=0$ is equivalent to
$$
J^{(i)}_{(\alpha)\beta j/\gamma}=0,\quad
J^{(i)}_{(\alpha)\beta j\vert k}=0,\quad
J^{(i)}_{(\alpha)\beta j}\vert^{(\gamma)}_{(k)}=0.
$$

In this context, we can prove the following

\begin{th}
The coefficients of an $h$-normal $\Gamma$-linear
connection $\nabla$ verify the identities
$$
\begin{array}{lll}\medskip
\bar G^\gamma_{\alpha\beta}=H^\gamma_{\alpha\beta},&\bar L^\alpha_{\beta j}=0,&
\bar C^{\alpha(\gamma)}_{\beta(j)}=0,\\\medskip
G^{(k)(\beta)}_{(\alpha)(i)\gamma}=\delta^\beta_\alpha G^k_{i\gamma}-\delta^k_i
H^\beta_{\alpha\gamma},&L^{(k)(\beta)}_{(\alpha)(i)j}=\delta^\beta_\alpha
L^k_{ij},&C^{(k)(\beta)(\gamma)}_{(\alpha)(i)(j)}=\delta^\beta_\alpha
C^{k(\gamma)}_{i(j)}.
\end{array}
$$
\end{th}
{\bf Proof.} The first three relations come from the definition of an $h$-normal
$\Gamma$-linear connection.

The condition $\nabla J=0$ implies locally that
$$
\left\{\begin{array}{l}\medskip
\displaystyle{h_{\beta\mu}G^{(i)(\mu)}_{(\alpha)(j)\gamma}=h_{\alpha\beta}G^i
_{j\gamma}+\delta^i_j\left[-{\partial h_{\alpha\beta}\over\partial t^\gamma}
+H_{\beta\gamma\alpha}\right]}\\\medskip
h_{\beta\mu}L^{(i)(\mu)}_{(\alpha)(j)}=h_{\alpha\beta}L^i_{jk}\\
h_{\beta\mu}C^{(i)(\mu)(\gamma)}_{(\alpha)(j)(k)}=h_{\alpha\beta}C^{i(\gamma)}
_{j(k)},
\end{array}\right.
$$
where $H_{\beta\gamma\alpha}=H^\mu_{\beta\gamma}h_{\mu\alpha}$ represent the
Christoffel symbols of the first kind attached to the pseudo-Riemannian metric
$h_{\alpha\beta}$. Contracting the above relations by $h^{\beta\varepsilon}$,
one obtains the last three identities of the theorem. \rule{5pt}{5pt}\medskip\\
\addtocounter{rem}{1}
{\bf Remarks \therem} i) The preceding theorem implies that an $h$-normal
$\Gamma$-linear on $E$ is determined just by four effective coefficients
\begin{equation}
\nabla\Gamma=(H^\gamma_{\alpha\beta},G^k_{i\gamma},L^k_{ij},C^{k(\gamma)}_
{i(j)}).
\end{equation}

ii) In the particular case $(T,h)=(R,\delta)$, a $\delta$-normal
$\Gamma$-linear connection identifies to the notion of $N$-linear connection
used in Lagrangian geometry \cite{5}.\medskip\\
\addtocounter{ex}{1}
{\bf Example \theex} The canonical Berwald $\Gamma_0$-linear connection associated
to the metrics pair $(h_{\alpha\beta},\varphi_{ij})$ is an $h$-normal
$\Gamma_0$-linear connection, defined by the local coefficients
$B\Gamma_0=(H^\gamma_{\alpha\beta},0,\gamma^k_{ij},0)$.

\section{Components of torsion and curvature d-tensors}

\setcounter{equation}{0}
\hspace{5mm} The study of the torsion {\bf T} and curvature {\bf R} d-tensor of an arbitrary
$\Gamma$-linear connection $\nabla$ on $E$ was made  in \cite{12}. In that
context, we proved that the torsion\linebreak d-tensor is
determined by twelve effective local torsion d-tensors, while the curvature
d-tensor of $\nabla$ is determined by eighteen local d-tensors.

Let us start with an $h$-normal $\Gamma$-linear connection
$\nabla$. Following the formulas described in \cite{12} and using the properties
of $\nabla$, it follows that the torsion d-tensors $\bar T^\mu_{\alpha\beta},\;\bar T^\mu_
{\alpha j}$ and $\bar P^{\mu(\beta)}_{\alpha(j)}$ vanish. Moreover, we
deduce that the following theorem holds good.
\begin{th}
The torsion d-tensor {\bf\em T} of an $h$-normal $\Gamma$-linear connection
$\nabla$, is determined by nine effective local d-tensors,
\begin{equation}
\begin{tabular}{|c|c|c|c|}
\hline
&$h_T$&$h_M$&$v$\\
\hline
$h_Th_T$&0&0&$R^{(m)}_{(\mu)\alpha\beta}$\\
\hline
$h_Mh_T$&0&$T^m_{\alpha j}$&$R^{(m)}_{(\mu)\alpha j}$\\
\hline
$h_Mh_M$&0&$T^m_{ij}$&$R^{(m)}_{(\mu)ij}$\\
\hline
$vh_T$&0&0&$P^{(m)\;\;(\beta)}_{(\mu)\alpha(j)}$\\
\hline
$vh_M$&0&$P^{m(\beta)}_{i(j)}$&$P^{(m)\;(\beta)}_{(\mu)i(j)}$\\
\hline
$vv$&0&0&$S^{(m)(\alpha)(\beta)}_{(\mu)(i)(j)}$\\
\hline
\end{tabular}
\end{equation}
where
$
\displaystyle{P^{(m)\;\;(\beta)}_{(\mu)\alpha(j)}={\partial M^{(m)}_{(\mu)
\alpha}\over\partial x^j_\beta}-\delta^\beta_\mu G^m_{j\alpha}+\delta^m_j
H^\beta_{\mu\alpha},}\quad
$
$\displaystyle{P^{(m)\;\;(\beta)}_{(\mu)i(j)}=
{\partial N^{(m)}_{(\mu)i}\over\partial x^j_\beta}-\delta^\beta_\mu L^m_{ji},}
$\medskip\linebreak
$\displaystyle{R^{(m)}_{(\mu)\alpha\beta}={\delta M^{(m)}_{(\mu)\alpha}\over
\delta t^\beta}-{\delta M^{(m)}_{(\mu)\beta}\over\delta t^\alpha},}
$
$\displaystyle{
R^{(m)}_{(\mu)\alpha j}={\delta M^{(m)}_{(\mu)\alpha}\over
\delta x^j}-{\delta N^{(m)}_{(\mu)j}\over\delta t^\alpha},}
$
$\displaystyle{
R^{(m)}_{(\mu)ij}={\delta N^{(m)}_{(\mu)i}\over\delta x^j}-{\delta N^{(m)}_
{(\mu)j}\over\delta x^i},}
$\medskip\\\\
$
S^{(m)(\alpha)(\beta)}_{(\mu)(i)(j)}=\delta^\alpha_\mu C^{m(\beta)}_{i(j)}-
\delta^\beta_\mu C^{m(\alpha)}_{j(i)}
$,
$
T^m_{\alpha j}=-G^m_{j\alpha}$,
$T^m_{ij}=L^m_{ij}-L^m_{ji}$,
$P^{m(\beta)}_{i(j)}=C^{m(\beta)}_{i(j)}$.
\end{th}
\addtocounter{ex}{1}
{\bf Remark \theex} For the Berwald $\Gamma_0$-linear connection associated to the
metrics $h_{\alpha\beta}$ and $\varphi_{ij}$, all torsion d-tensors vanish,
except
\begin{equation}
R^{(m)}_{(\mu)\alpha\beta}=-H^\gamma_{\mu\alpha\beta}x^m_\gamma,\quad
R^{(m)}_{(\mu)ij}=r^m_{ijl}x^l_\mu,
\end{equation}
where $H^\gamma_{\mu\alpha\beta}$ (resp. $r^m_{ijl}$) are the curvature tensors
of the metric $h_{\alpha\beta}$ (resp. $\varphi_{ij}$).\medskip

The form of expressions of local curvature d-tensors from the general case of
a $\Gamma$-linear connection \cite {12}, and again the properties of the $h$-normal
$\Gamma$-linear connection $\nabla$, imply a reduction (from eighteen to seven)
of the number of the effective curvature d-tensors attached to an $h$-normal $\Gamma$-linear
connection. Consequently, we obtain
\begin{th}
The curvature d-tensor {\bf\em R} of an $h$-normal $\Gamma$-linear connection
$\nabla$, is characterized by seven effective local d-tensors,
\begin{equation}
\begin{tabular}{|c|c|c|c|}
\hline
&$h_T$&$h_M$&$v$\\
\hline
$h_Th_T$&$H^\alpha_{\eta\beta\gamma}$&$R^l_{i\beta\gamma}$&
$R^{(l)(\alpha)}_{(\eta)(i)\beta\gamma}=\delta^\alpha_\eta R^l_{i\beta\gamma}
+\delta^l_iH^\alpha_{\eta\beta\gamma}$\\
\hline
$h_Mh_T$&0&$R^l_{i\beta k}$&$R^{(l)(\alpha)}_{(\eta)(i)\beta k}=\delta^\alpha
_\eta R^l_{i\beta k}$\\
\hline
$h_Mh_M$&0&$R^l_{ijk}$&$R^{(l)(\alpha)}_{(\eta)(i)jk}=\delta^\alpha_\eta R^l_{
ijk}$\\
\hline
$vh_T$&0&$P^{l\;\;(\gamma)}_{i\beta(k)}$&$P^{(l)(\alpha)\;\;(\gamma)}_{(\eta)
(i)\beta(k)}=\delta^\alpha_\eta P^{l\;\;(\gamma)}_{i\beta(k)}$\\
\hline
$vh_M$&0&$P^{l\;(\gamma)}_{ij(k)}$&$P^{(l)(\alpha)\;(\gamma)}_{(\eta)
(i)j(k)}=\delta^\alpha_\eta P^{l\;(\gamma)}_{ij(k)}$\\
\hline
$vv$&0&$S^{l(\beta)(\gamma)}_{i(j)(k)}$&$S^{(l)(\alpha)(\beta)(\gamma)}_
{(\eta)(i)(j)(k)}=\delta^\alpha_\eta S^{l(\beta)(\gamma)}_{i(j)(k)}$\\
\hline
\end{tabular}
\end{equation}
where\medskip

$
\displaystyle{H^\alpha_{\eta\beta\gamma}={\partial H^\alpha_{\eta\beta}\over
\partial t^\gamma}-{\partial H^\alpha_{\eta\gamma}\over\partial t^\beta}+
H^\mu_{\eta\beta}H^\alpha_{\mu\gamma}-H^\mu_{\eta\gamma}H^\alpha_{\mu\beta},}
$\medskip

$
\displaystyle{R^l_{i\beta\gamma}={\delta G^l_{i\beta}\over\delta t^\gamma}-
{\delta G^l_{i\gamma}\over\delta t^\beta}+G^m_{i\beta}G^l_{m\gamma}-
G^m_{i\gamma}G^l_{m\beta}+C^{l(\mu)}_{i(m)}R^{(m)}_{(\mu)\beta\gamma},}
$\medskip

$
\displaystyle{R^l_{i\beta k}={\delta G^l_{i\beta}\over\delta x^k}-
{\delta L^l_{ik}\over\delta t^\beta}+G^m_{i\beta}L^l_{mk}-
L^m_{ik}G^l_{m\beta}+C^{l(\mu)}_{i(m)}R^{(m)}_{(\mu)\beta k},}
$\medskip

$
\displaystyle{R^l_{ijk}={\delta L^l_{ij}\over\delta x^k}-
{\delta L^l_{ik}\over\delta x^j}+L^m_{ij}L^l_{mk}-L^m_{ik}L^l_{mj}+
C^{l(\mu)}_{i(m)}R^{(m)}_{(\mu)jk},}
$\medskip

$
\displaystyle{P^{l\;\;(\gamma)}_{i\beta(k)}={\partial G^l_{i\beta}\over\partial
x^k_\gamma}-C^{l(\gamma)}_{i(k)/\beta}+C^{l(\mu)}_{i(m)}P^{(m)\;\;(\gamma)}_
{(\mu)\beta(k)},}
$\medskip

$
\displaystyle{P^{l\;(\gamma)}_{ij(k)}={\partial L^l_{ij}\over\partial
x^k_\gamma}-C^{l(\gamma)}_{i(k)\vert j}+C^{l(\mu)}_{i(m)}P^{(m)\;(\gamma)}_
{(\mu)j(k)},}
$\medskip

$
\displaystyle{S^{l(\beta)(\gamma)}_{i(j)(k)}={\partial C^{l(\beta)}_{i(j)}
\over\partial x^k_\gamma}-{\partial C^{l(\gamma)}_{i(k)}\over\partial x^j_
\beta}+C^{m(\beta)}_{i(j)}C^{l(\gamma)}_{m(k)}-
C^{m(\gamma)}_{i(k)}C^{l(\beta)}_{m(j)}.}
$
\end{th}
\addtocounter{rem}{1}
{\bf Remark \therem} In the case of the Berwald $\Gamma_0$-linear connection
associated to the metric pair $(h_{\alpha\beta},\varphi_{ij})$, all curvature
d-tensors vanish, except $H^\delta_{\alpha\beta\gamma}$ and $R^l_{ijk}=r^l_
{ijk}$, where $r^l_{ijk}$ are the curvature tensors of the metric $\varphi_{ij}$.

\section{Ricci identities. Deflection d-tensors identities}

\setcounter{equation}{0}
\hspace{5mm} The Ricci identities of a $\Gamma$-linear  connection are described
in  \cite{12}. In the particular case of an $h$-normal
$\Gamma$-linear connection, these simplify because the number and
the form of the torsion and curvature d-tensors reduced. A meaningful reduction
of these identities can be obtained, considering the more particular case of an
$h$-normal \linebreak $\Gamma$-linear connection $\nabla$ of {\it Cartan type},
(i. e. , $L^i_{jk}=L^i_{kj}$ and $C^{i(\gamma)}_{j(k)}=C^{i(\gamma)}_{k(j)}$).
In that case, the condition $L^i_{jk}=L^i_{kj}$ implies $T^i_{jk}=0$.
Consequently, we have
\begin{th}
The following Ricci identities of an $h$-normal $\Gamma$-linear connection of
Cartan type, are true:\medskip\\
$(h_T)\mbox{\hspace{6mm}}
\left\{\begin{array}{l}\medskip
X^\alpha_{/\beta/\gamma}-X^\alpha_{/\gamma/\beta}=X^\mu H^\alpha_{\mu\beta
\gamma}-X^\alpha\vert^{(\mu)}_{(m)}R^{(m)}_{(\mu)\beta\gamma}\\\medskip
X^\alpha_{/\beta\vert k}-X^\alpha_{\vert k/\beta}=-X^\alpha_{\vert m}
T^m_{\beta k}-X^\alpha\vert^{(\mu)}_{(m)}R^{(m)}_{(\mu)\beta k}\\\medskip
X^\alpha_{\vert j\vert k}-X^\alpha_{\vert k\vert j}=-X^\alpha\vert^{(\mu)}_
{(m)}R^{(m)}_{(\mu)jk}\\\medskip
X^\alpha_{/\beta}\vert^{(\gamma)}_{(k)}-X^\alpha\vert^{(\gamma)}_{(k)/\beta}=
-X^\alpha\vert^{(\mu)}_{(m)}P^{(m)\;\;(\gamma)}_{(\mu)\beta(k)}\\\medskip
X^\alpha_{\vert j}\vert^{(\gamma)}_{(k)}-X^\alpha\vert^{(\gamma)}_{(k)\vert j}
=-X^\alpha_{\vert m}C^{m(\gamma)}_{j(k)}-X^\alpha\vert^{(\mu)}_{(m)}P^{(m)\;
(\gamma)}_{(\mu)j(k)}\\\medskip
X^\alpha\vert^{(\beta)}_{(j)}\vert^{(\gamma)}_{(k)}-X^\alpha\vert^{(\gamma)}_
{(k)}\vert^{(\beta)}_{(j)}=-X^\alpha\vert^{(\mu)}_{(m)}S^{(m)(\beta)(\gamma)}_
{(\mu)(j)(k)},
\end{array}\right.
$\\\medskip
$(h_M)\mbox{\hspace{5mm}}
\left\{\begin{array}{l}\medskip
X^i_{/\beta/\gamma}-X^i_{/\gamma/\beta}=X^m R^i_{m\beta\gamma}-X^i\vert^
{(\mu)}_{(m)}R^{(m)}_{(\mu)\beta\gamma}\\\medskip
X^i_{/\beta\vert k}-X^i_{\vert k/\beta}=X^mR^i_{m\beta k}-X^i_{\vert m}
T^m_{\beta k}-X^i\vert^{(\mu)}_{(m)}R^{(m)}_{(\mu)\beta k}\\\medskip
X^i_{\vert j\vert k}-X^i_{\vert k\vert j}=X^mR^i_{mjk}-X^i\vert^{(\mu)}_
{(m)}R^{(m)}_{(\mu)jk}\\\medskip
X^i_{/\beta}\vert^{(\gamma)}_{(k)}-X^i\vert^{(\gamma)}_{(k)/\beta}=
X^mP^{i\;\;(\gamma)}_{m\beta(k)}-X^i\vert^{(\mu)}_{(m)}P^{(m)\;\;(\gamma)}_
{(\mu)\beta(k)}\\\medskip
X^i_{\vert j}\vert^{(\gamma)}_{(k)}-X^i\vert^{(\gamma)}_{(k)\vert j}=
X^mP^{i\;(\gamma)}_{mj(k)}-X^i_{\vert m}C^{m(\gamma)}_{j(k)}-X^i\vert^{(\mu)}
_{(m)}P^{(m)\;(\gamma)}_{(\mu)j(k)}\\\medskip
X^i\vert^{(\beta)}_{(j)}\vert^{(\gamma)}_{(k)}-X^i\vert^{(\gamma)}_{(k)}
\vert^{(\beta)}_{(j)}=X^mS^{i(\beta)(\gamma)}_{m(j)(k)}-X^i\vert^{(\mu)}_
{(m)}S^{(m)(\beta)(\gamma)}_{(\mu)(j)(k)},
\end{array}\right.
$\\\medskip
$(v)\mbox{\hspace{8mm}}
\left\{\begin{array}{l}\medskip
X^{(i)}_{(\alpha)/\beta/\gamma}-X^{(i)}_{(\alpha)/\gamma/\beta}=X^{(m)}_
{(\alpha)}R^i_{m\beta\gamma}-X^{(i)}_{(\mu)}H^\mu_{\alpha\beta\gamma}-
X^{(i)}_{(\alpha)}\vert^{(\mu)}_{(m)}R^{(m)}_{(\mu)\beta\gamma}\\\medskip
X^{(i)}_{(\alpha)/\beta\vert k}-X^{(i)}_{(\alpha)\vert k/\beta}=X^{(m)}_
{(\alpha)}R^i_{m\beta k}-X^{(i)}_{(\alpha)\vert m}T^m_{\beta k}-X^{(i)}_
{(\alpha)}\vert^{(\mu)}_{(m)}R^{(m)}_{(\mu)\beta k}\\\medskip
X^{(i)}_{(\alpha)\vert j\vert k}-X^{(i)}_{(\alpha)\vert k\vert j}=X^{(m)}_
{(\alpha)}R^i_{mjk}-X^{(i)}_{(\alpha)}\vert^{(\mu)}_{(m)}R^{(m)}_{(\mu)jk}
\\\medskip
X^{(i)}_{(\alpha)/\beta}\vert^{(\gamma)}_{(k)}-X^{(i)}_{(\alpha)}\vert^
{(\gamma)}_{(k)/\beta}=X^{(m)}_{(\alpha)}P^{i\;\;(\gamma)}_{m\beta(k)}-X^{(i)}
_{(\alpha)}\vert^{(\mu)}_{(m)}P^{(m)\;\;(\gamma)}_{(\mu)\beta(k)}\\\medskip
X^{(i)}_{(\alpha)\vert j}\vert^{(\gamma)}_{(k)}-X^{(i)}_{(\alpha)}\vert^
{(\gamma)}_{(k)\vert j}=X^{(m)}_{(\alpha)}P^{i\;(\gamma)}_{mj(k)}-X^{(i)}_
{(\alpha)\vert m}C^{m(\gamma)}_{j(k)}-X^{(i)}_{(\alpha)}\vert^{(\mu)}_{(m)}
P^{(m)\;(\gamma)}_{(\mu)j(k)}\\\medskip
X^{(i)}_{(\alpha)}\vert^{(\beta)}_{(j)}\vert^{(\gamma)}_{(k)}-X^{(i)}_{(\alpha)}
\vert^{(\gamma)}_{(k)}\vert^{(\beta)}_{(j)}=X^{(m)}_{(\alpha)}S^{i(\beta)
(\gamma)}_{m(j)(k)}-X^{(i)}_{(\alpha)}\vert^{(\mu)}_{(m)}S^{(m)(\beta)
(\gamma)}_{(\mu)(j)(k)},
\end{array}\right.
$\\
where $\displaystyle{X=X^\alpha{\delta\over\delta t^\alpha}+X^i{\delta\over
\delta x^i}+X^{(i)}_{(\alpha)}{\partial\over\partial x^i_\alpha}}$ is an
arbitrary d-vector field on $J^1(T,M)$.
\end{th}

In the sequel, let us consider the canonical Liouville d-tensor {\bf C}$
\displaystyle{=x^i_\alpha{\partial\over\partial x^i_\alpha}}$. We construct
the {\it deflection d-tensors associated to the $h$-normal $\Gamma$-linear
connection $\nabla$}, setting
\begin{equation}
\bar D^{(i)}_{(\alpha)\beta}=x^i_{\alpha/\beta},\quad
D^{(i)}_{(\alpha)k}=x^i_{\alpha\vert k},\quad
d^{(i)(\gamma)}_{(\alpha)(k)}=x^i_\alpha\vert^{(\gamma)}_{(k)},
\end{equation}
where $"_{/\beta}"$, $"_{\vert k}"$ and $"\vert^{(\gamma)}_{(k)}$ are the local
covariant derivatives induced by $\nabla$.
By a direct calculation, the deflection d-tensors get the expressions
\begin{equation}
\left\{\begin{array}{l}\medskip
\bar D^{(i)}_{(\alpha)\beta}=-M^{(i)}_{(\alpha)\beta}+G^i_{m\beta}x^m_\alpha-
H^\mu_{\alpha\beta}x^i_\mu\\\medskip
D^{(i)}_{(\alpha)j}=-N^{(i)}_{(\alpha)j}+L^i_{mj}x^m_\alpha\\\medskip
d^{(i)(\beta)}_{(\alpha)(j)}=\delta^i_j\delta^\beta_\alpha+C^{i(\beta)}_
{m(j)}x^m_\alpha.
\end{array}\right.
\end{equation}

Applying the $(v)$-set of  the Ricci identities to the components of the Liouville
vector field, we obtain
\begin{th}
The deflection d-tensors,  attached to the $h$-normal $\Gamma$-linear connection
$\nabla$, verify the identities:
\begin{equation}
\left\{\begin{array}{l}\medskip
\bar D^{(i)}_{(\alpha)\beta/\gamma}-\bar D^{(i)}_{(\alpha)\gamma/\beta}=
x^m_\alpha R^i_{m\beta\gamma}-x^i_\mu H^\mu_{\alpha\beta\gamma}-
d^{(i)(\mu)}_{(\alpha)(m)}R^{(m)}_{(\mu)\beta\gamma}\\\medskip
\bar D^{(i)}_{(\alpha)\beta\vert k}-D^{(i)}_{(\alpha)k/\beta}=x^m_\alpha
R^i_{m\beta k}-D^{(i)}_{(\alpha)m}T^m_{\beta k}-d^{(i)(\mu)}_{(\alpha)(m)}
R^{(m)}_{(\mu)\beta k}\\\medskip
D^{(i)}_{(\alpha)j\vert k}-D^{(i)}_{(\alpha)k\vert j}=x^m_\alpha R^i_{mjk}-
d^{(i)(\mu)}_{(\alpha)(m)}R^{(m)}_{(\mu)jk}
\\\medskip
\bar D^{(i)}_{(\alpha)\beta}\vert^{(\gamma)}_{(k)}-d^{(i)(\gamma)}_{(\alpha)
(k)/\beta}=x^m_\alpha P^{i\;\;(\gamma)}_{m\beta(k)}-d^{(i)(\mu)}
_{(\alpha)(m)}P^{(m)\;\;(\gamma)}_{(\mu)\beta(k)}\\\medskip
D^{(i)}_{(\alpha)j}\vert^{(\gamma)}_{(k)}-d^{(i)(\gamma)}_{(\alpha)
(k)\vert j}=x^m_\alpha P^{i\;(\gamma)}_{mj(k)}-D^{(i)}_{(\alpha)m}
C^{m(\gamma)}_{j(k)}-d^{(i)(\mu)}_{(\alpha)(m)}P^{(m)\;(\gamma)}_
{(\mu)j(k)}\\\medskip
d^{(i)(\beta)}_{(\alpha)(j)}\vert^{(\gamma)}_{(k)}-d^{(i)(\gamma)}_{(\alpha)
(k)}\vert^{(\beta)}_{(j)}=x^m_\alpha S^{i(\beta)(\gamma)}_{m(j)(k)}-
d^{(i)(\mu)}_{(\alpha)(m)}S^{(m)(\beta)(\gamma)}_{(\mu)(j)(k)}.
\end{array}\right.
\end{equation}
\end{th}
\addtocounter{rem}{1}
{\bf Remark \therem} The importance of the deflection d-tensors identities is
emphasized in \cite{7}, \cite{9}, where is developped the {\it (generalized)
metrical multi-time Lagrangian}  geometry of physical fields  on $J^1(T,M)$.
In that context, the deflection d-tensors identities  are used in the
description of the  Maxwell equations which govern the electromagnetic field
of a {\it (generalized) metrical multi-time Lagrange space} \cite{7}, \cite{11}.

\section{Bianchi identities}

\setcounter{equation}{0}
\hspace{5mm} From the general theory of linear connections  on a vector bundle
$E$, is known that the torsions {\bf T} and the curvature {\bf R} of a linear
connection $\nabla$ are not independent. They verify the Bianchi identities,
whose expressions, in a local basis $(X_A)$ of ${\cal X}(E)$, are \cite{5}, \cite{12}
\begin{equation}
\left\{\begin{array}{l}\medskip
\sum_{\{A,B,C\}}\{R^F_{ABC}-T^F_{AB:C}-T^G_{AB}T^F_{CG}\}=0\\
\sum_{\{A,B,C\}}\{R^F_{DAB:C}+T^G_{AB}R^F_{DAG}\}=0,
\end{array}\right.
\end{equation}
where {\bf R}$(X_A,X_B)X_C=R^D_{CBA}X_D$, {\bf T}$(X_A,X_B)=T^D_{BA}X_D$ and
"$_{:C}$" represents the local covariant derivative induced  by $\nabla$.

In our context, we have $E=J^1(T,M)$. Let $\Gamma=(M^{(i)}_{(\alpha)\beta},
N^{(i)}_{(\alpha)j})$ be a fixed nonlinear  connection on $E$, and
$\displaystyle{(X_A)=\left({\delta\over\delta t^\alpha},{\delta\over\delta
x^i},{\partial\over\partial x^i_\alpha}\right)}$ its adapted basis.
In the sequel, we try to rewrite the above Bianchi identities for
an $h$-normal $\Gamma$-linear connection $\nabla$ of Cartan type on $E$. In this sense,
taking into account that the indices $A,B,\ldots$ are of type $\left\{\alpha,
i,{(\alpha)\atop(i)}\right\}$, it follows that the covariant derivative $"_{:A}"$
becomes one of the coveriant derivatives  $"_{/\alpha}"$, $"_{\vert i}"$
or $"\vert^{(\alpha)}_{(i)}"$. Consequently, we deduce
\begin{th}
The following thirty effective Bianchi identities of the $h$-normal\linebreak
$\Gamma$-linear connection $\nabla$ of Cartan type are true:
\medskip\\
$(1)\mbox{\hspace{3mm}}
\left\{\begin{array}{ll}\medskip
1.\;1&\sum_{\{\alpha,\beta,\gamma\}}H^\delta_{\alpha\beta\gamma}=0
\\\medskip
1.\;2&{\cal A}_{\{\alpha,\beta\}}\left\{T^l_{\alpha m}T^m_{\beta k}-T^l_{\alpha k/
\beta}\right\}=R^l_{k\alpha\beta}-C^{l(\mu)}_{k(m)}R^{(m)}_{(\mu)\alpha\beta}
\\\medskip
1.\;3&{\cal A}_{\{j,k\}}\left\{C^{l(\mu)}_{k(m)}R^{(m)}_{(\mu)\alpha j}+R^l_
{j\alpha k}+T^l_{\alpha j\vert k}\right\}=0
\\\medskip
1.\;4&\sum_{\{i,j,k\}}\left\{C^{l(\mu)}_{k(m)}R^{(m)}_{(\mu)ij}-R^l_{ijk}\right\}=0,
\end{array}\right.
$\medskip\linebreak
$(2)\mbox{\hspace{3mm}}
\left\{\begin{array}{ll}\medskip
2.\;1&\sum_{\{\alpha,\beta,\gamma\}}\left\{R^{(l)}_{(\delta)\alpha\beta/\gamma}+
P^{(l)\;\;(\mu)}_{(\delta)\gamma(m)}R^{(m)}_{(\mu)\alpha\beta}\right\}=0
\\\medskip
2.\;2&{\cal A}_{\{\alpha,\beta\}}\left\{R^{(l)}_{(\delta)\alpha k/\beta}+P^{(l)\;\;
(\mu)}_{(\delta)\beta(m)}R^{(m)}_{(\mu)\alpha k}+R^{(l)}_
{(\delta)\beta m}T^m_{\alpha k}\right\}=R^{(l)}_{(\delta)\alpha\beta\vert k}+\\\medskip
&\mbox{\hspace{75mm}}+P^{(l)\;\;(\mu)}_{(\delta)k(m)}R^{(m)}_{(\mu)\alpha\beta}
\\\medskip
2.\;3&{\cal A}_{\{j,k\}}\left\{R^{(l)}_{(\delta)\alpha j\vert k}+P^{(l)\;(\mu)}
_{(\delta)k(m)}R^{(m)}_{(\mu)\alpha j}+R^{(l)}_{(\delta)km}T^m_{\alpha j}\right\}=
-R^{(l)}_{(\delta)\alpha j\vert k}-\\\medskip
&\mbox{\hspace{72mm}}-P^{(l)\;\;(\mu)}_{(\delta)\alpha(m)}R^{(m)}_{(\mu)jk}
\\\medskip
2.\;4&\sum_{\{i,j,k\}}\left\{R^{(l)}_{(\delta)ij\vert k}+P^{(l)\;(\mu)}_{(\delta)k
(m)}R^{(m)}_{(\mu)ij}\right\}=0,
\end{array}\right.
$\medskip\linebreak
$(3)\mbox{\hspace{3mm}}
\left\{\begin{array}{ll}\medskip
3.\;1&T^l_{\alpha k}\vert^{(\varepsilon)}_{(p)}-C^{l(\varepsilon)}_{m(p)}
T^m_{\alpha k}+P^{l\;\;(\varepsilon)}_{k\alpha(p)}-C^{l(\varepsilon)}_
{k(p)/\alpha}-C^{l(\mu)}_{k(m)}P^{(m)\;\;(\varepsilon)}_{(\mu)\alpha(p)}=0
\\\medskip
3.\;2&{\cal A}_{\{j,k\}}\left\{C^{l(\varepsilon)}_{j(p)\vert k}+C^{l(\mu)}_{k(m)}
P^{(m)\;(\varepsilon)}_{(\mu)j(p)}+P^{l\;(\varepsilon)}_{jk(p)}\right\}=0,
\end{array}\right.
$\medskip\linebreak
$(4)\mbox{\hspace{3mm}}
\left\{\begin{array}{ll}\medskip
4.\;1&{\cal A}_{\{\alpha,\beta\}}\left\{P^{(l)\;\;(\varepsilon)}_{(\delta)
\alpha(p)/\beta}+P^{(l)\;\;(\mu)}_{(\delta)\beta(m)}P^{(m)\;\;(\varepsilon)}_
{(\mu)\alpha(p)}\right\}=R^{(l)}_{(\delta)\alpha\beta}\vert^{(\varepsilon)}_
{(p)}-\\\medskip
&\mbox{\hspace{55mm}}-R^{(l)(\varepsilon)}_{(\delta)(p)\alpha\beta}+S^{(l)(\varepsilon)(\mu)}_
{(\delta)(p)(m)}R^{(m)}_{(\mu)\alpha\beta}
\\\medskip
4.\;2&{\cal A}_{\{\alpha,k\}}\left\{P^{(l)\;\;(\varepsilon)}_{(\delta)
\alpha(p)\vert k}+P^{(l)\;\;(\mu)}_{(\delta)k(m)}P^{(m)\;\;(\varepsilon)}_
{(\mu)\alpha(p)}\right\}=R^{(l)}_{(\delta)\alpha k}\vert^{(\varepsilon)}_
{(p)}-\\\medskip
&\mbox{\hspace{10mm}}-R^{(l)(\varepsilon)}_{(\delta)(p)\alpha k}+S^{(l)(\varepsilon)(\mu)}_
{(\delta)(p)(m)}R^{(m)}_{(\mu)\alpha k}+R^{(l)}_{(\delta)\alpha m}C^{m(\varepsilon)}_
{k(p)}-T^m_{\alpha k}P^{(l)\;(\varepsilon)}_{(\delta)m(p)}
\\\medskip
4.\;3&{\cal A}_{\{j,k\}}\left\{P^{(l)\;\;(\varepsilon)}_{(\delta)j(p)\vert k}
+P^{(l)\;\;(\mu)}_{(\delta)k(m)}P^{(m)\;\;(\varepsilon)}_{(\mu)j(p)}+R^{(l)}_
{(\delta)km}C^{m(\varepsilon)}_{j(p)}\right\}=\\\medskip
&\mbox{\hspace{35mm}}=R^{(l)}_{(\delta)jk}\vert^{(\varepsilon)}_{(p)}-R^{(l)(\varepsilon)}_
{(\delta)(p)jk}+S^{(l)(\varepsilon)(\mu)}_{(\delta)(p)(m)}R^{(m)}_{(\mu)jk},
\end{array}\right.
$\medskip\linebreak
$(5)\mbox{\hspace{3mm}}
\left\{\begin{array}{ll}\medskip
5.\;1&{\cal A}_{\left\{{(\beta)\atop(j)},{(\gamma)\atop(k)}\right\}}\left\{
C^{l(\beta)}_{i(j)}\vert^{(\gamma)}_{(k)}+C^{m(\gamma)}_{i(k)}C^{l(\beta)}_
{m(j)}\right\}=S^{l(\beta)(\gamma)}_{i(j)(k)}-C^{l(\mu)}_{i(m)}S^{(m)(\beta)
(\gamma)}_{(\mu)(j)(k)},
\end{array}\right.
$\medskip\linebreak
$(6)\mbox{\hspace{3mm}}
\left\{\begin{array}{ll}\medskip
6.\;1&{\cal A}_{\left\{{(\beta)\atop(j)},{(\gamma)\atop(k)}\right\}}\left\{
P^{(l)\;\;(\beta)}_{(\delta)\alpha(j)}\vert^{(\gamma)}_{(k)}+P^{(m)\;\;(\beta)}
_{(\mu)\alpha(j)}S^{(l)(\gamma)(\mu)}_{(\delta)(k)(m)}+P^{(l)(\beta)\;\;(\gamma)}_
{(\delta)(j)\alpha(k)}\right\}=\\\medskip
&\hbox{\hspace{45mm}}=-S^{(l)(\beta)(\gamma)}_{(\delta)(j)(k)/\alpha}-
S^{(m)(\beta)(\gamma)}_{(\mu)(j)(k)}P^{(l)\;\;(\mu)}_{(\delta)\alpha(m)}
\\\medskip
6.\;2&{\cal A}_{\left\{{(\beta)\atop(j)},{(\gamma)\atop(k)}\right\}}\left\{
P^{(l)\;(\beta)}_{(\delta)i(j)}\vert^{(\gamma)}_{(k)}+P^{(m)\;(\beta)}
_{(\mu)i(j)}S^{(l)(\gamma)(\mu)}_{(\delta)(k)(m)}+P^{(l)(\beta)\;(\gamma)}_
{(\delta)(j)i(k)}\right\}=\\\medskip
&\hbox{\hspace{45mm}}=-S^{(l)(\beta)(\gamma)}_{(\delta)(j)(k)\vert i}-
S^{(m)(\beta)(\gamma)}_{(\mu)(j)(k)}P^{(l)\;(\mu)}_{(\delta)i(m)},
\end{array}\right.
$\medskip\linebreak
$(7)\mbox{\hspace{3mm}}
\left\{\begin{array}{ll}\medskip
7.\;1&
\sum_{\left\{{(\alpha)\atop(i)},{(\beta)\atop(j)},{(\gamma)\atop(k)}\right\}}
\left\{S^{(l)(\alpha)(\beta)}_{(\delta)(i)(j)}\vert^{(\gamma)}_{(k)}+
S^{(m)(\alpha)(\beta)}_{(\mu)(i)(j)}S^{(l)(\gamma)(\mu)}_{(\delta)(k)(m)}-
S^{(l)(\alpha)(\beta)(\gamma)}_{(\delta)(i)(j)(k)}\right\}=0,
\end{array}\right.
$\medskip\linebreak
$(8)\mbox{\hspace{3mm}}
\left\{\begin{array}{ll}\medskip
8.\;1&\sum_{\{\alpha,\beta,\gamma\}}H^\delta_{\varepsilon\alpha\beta/\gamma}=0
\\\medskip
8.\;2&H^\delta_{\varepsilon\alpha\beta\vert k}=0
\\\medskip
8.\;3&\sum_{\{i,j,k\}}R^{(m)}_{(\mu)ij}P^{(\delta)\;(\mu)}_{(\varepsilon)k(m)}=0
\\\medskip
8.\;4&\sum_{\{\alpha,\beta,\gamma\}}\left\{R^l_{p\alpha\beta/\gamma}-R^{(m)}_
{(\mu)\alpha\beta}P^{l\;\;(\mu)}_{p\gamma(m)}\right\}=0
\\\medskip
8.\;5&{\cal A}_{\{\alpha,\beta\}}\left\{R^l_{p\alpha k/\beta}+R^{(m)}_{(\mu)
\alpha k}P^{l\;\;(\mu)}_{p\beta(m)}-T^m_{\alpha k}R^l_{p\beta m}\right\}=
R^l_{p\alpha\beta}+R^{(m)}_{(\mu)\alpha\beta}P^{l\;(\mu)}_{pk(m)}
\\\medskip
8.\;6&{\cal A}_{\{j,k\}}\left\{R^l_{p\alpha j\vert k}+R^{(m)}_{(\mu)\alpha j}
P^{l\;\;(\mu)}_{pk(m)}-T^m_{\alpha j}R^l_{pkm}\right\}=
-R^l_{pjk/\alpha}+R^{(m)}_{(\mu)\alpha k}P^{l\;(\mu)}_{pj(m)}
\\\medskip
8.\;7&\sum_{\{i,j,k\}}\left\{R^l_{pij\vert k}-R^{(m)}_{(\mu)ij}P^{l\;(\mu)}_
{pk(m)}\right\}=0,
\end{array}\right.
$\medskip\linebreak
$(9)\mbox{\hspace{3mm}}
\left\{\begin{array}{ll}\medskip
9.\;1&{\cal A}_{\{\alpha,\beta\}}\left\{P^{l\;\;(\varepsilon)}_{i\alpha(p)/\beta}-
P^{(m)\;\;(\varepsilon)}_{(\mu)\alpha(p)}P^{l\;\;(\mu)}_{i\beta(m)}\right\}=
R^l_{i\alpha\beta}\vert^{(\varepsilon)}_{(p)}+R^{(m)}_{(\mu)\alpha\beta}
S^{l(\varepsilon)(\mu)}_{i(p)(m)}
\\\medskip
9.\;2&{\cal A}_{\{\alpha,k\}}\left\{P^{l\;\;(\varepsilon)}_{i\alpha(p)\vert k}-
P^{(m)\;\;(\varepsilon)}_{(\mu)\alpha(p)}P^{l\;\;(\mu)}_{ik(m)}\right\}=
R^l_{i\alpha k}\vert^{(\varepsilon)}_{(p)}-R^{(m)}_{(\mu)\alpha k}
S^{l(\varepsilon)(\mu)}_{i(p)(m)}-\\\medskip
&\mbox{\hspace{60mm}}-C^{m(\varepsilon)}_{k(p)}R^l_{i\alpha  m}+T^m_{\alpha k}
P^{l\;\;(\varepsilon)}_{im(p)}
\\\medskip
9.\;3&{\cal A}_{\{j,k\}}\left\{P^{l\;\;(\varepsilon)}_{ij(p)\vert k}-P^{(m)\;\;
(\varepsilon)}_{(\mu)j(p)}P^{l\;\;(\mu)}_{ik(m)}-C^{m(\varepsilon)}_{j(p)}R^l
_{ikm}\right\}=R^l_{ijk}\vert^{(\varepsilon)}_{(p)}+\\\medskip
&\mbox{\hspace{74mm}}+R^{(m)}_{(\mu)jk}S^{l(\varepsilon)(\mu)}_{i(p)(m)},
\end{array}\right.
$\medskip\linebreak
$(10)\mbox{\hspace{1mm}}
\left\{\begin{array}{ll}\medskip
10.\;1&{\cal A}_{\left\{{(\beta)\atop(j)},{(\gamma)\atop(k)}\right\}}\left\{
P^{l\;\;(\beta)}_{p\alpha(j)}\vert^{(\gamma)}_{(k)}-P^{(m)\;\;(\beta)}_{(\mu)
\alpha(j)}S^{l(\gamma)(\mu)}_{p(k)(m)}\right\}=S^{l(\beta)(\gamma)}_{p(j)(k)/
\alpha}+\\\medskip
&\mbox{\hspace{66mm}}+S^{(m)(\beta)(\gamma)}_{(\mu)(j)(k)}P^{l\;\;(\mu)}_{p\alpha(m)}
\\\medskip
10.\;2&{\cal A}_{\left\{{(\beta)\atop(j)},{(\gamma)\atop(k)}\right\}}\left\{
P^{l\;\;(\beta)}_{pi(j)}\vert^{(\gamma)}_{(k)}-P^{(m)\;\;(\beta)}_{(\mu)
i(j)}S^{l(\gamma)(\mu)}_{p(k)(m)}+C^{m(\beta)}_{i(j)}P^{l\;\;(\gamma)}_
{pm(k)}\right\}=\\\medskip
&\mbox{\hspace{50mm}}=S^{l(\beta)(\gamma)}_{p(j)(k)\vert i}+
S^{(m)(\beta)(\gamma)}_{(\mu)(j)(k)}P^{l\;\;(\mu)}_{pi(m)},
\end{array}\right.
$\medskip\linebreak
$(11)\mbox{\hspace{1mm}}
\left\{\begin{array}{ll}\medskip
11.\;1&
\sum_{\left\{{(\alpha)\atop(i)},{(\beta)\atop(j)},{(\gamma)\atop(k)}\right\}}
\left\{S^{l(\alpha)(\beta)}_{p(i)(j)}\vert^{(\gamma)}_{(k)}+
S^{(m)(\alpha)(\beta)}_{(\mu)(i)(j)}S^{l(\gamma)(\mu)}_{p(k)(m)}\right\}=0,
\end{array}\right.
$\medskip\linebreak
where, if $\{A,B,C\}$ are indices of type $\left\{\alpha,i,{(\alpha)\atop(i)}
\right\}$ then $\sum_{\{A,B,C\}}$ means the cyclic sum and ${\cal A}_{\{A,B\}}$
means alternate sum.
\end{th}
\addtocounter{rem}{1}
{\bf Remarks \therem} i) In the particular case $(T,h)=(R,\delta)$, the last
identity of every  above set of the Bianchi identities reduces to the
one of the classical eleven Bianchi identities of an $N$-linear connection
from the Lagrange geometry \cite{5}.

ii) The Bianchi identities of an $h$-normal $\Gamma$-linear connection of Cartan
type are used  in the description of the Maxwell equations and the conservation
laws of the Einstein equations of the gravitational potentials from the
background of the (generalized) metrical multi-time Lagrange geometry of
physical fields \cite{7}, \cite{9}.\medskip\\
{\bf Acknoledgements.} The author would like to thank to the reviewers of
Journal of the Mathematical Society of Japan for their valuable comments
upon a previous version of this paper.

\begin{center}
University POLITEHNICA of Bucharest\\
Department of Mathematics I\\
Splaiul Independentei 313\\
77206 Bucharest, Romania\\
e-mail: mircea@mathem.pub.ro\\
Fax: (401)411.53.65.
\end{center}

\end{document}